\documentclass[twocolumn]{article}
\usepackage{natbib}

\catcode`@=11

\def\intro{\@startsection {section}{10}{\z@}{-3.5ex plus -1ex minus
 -.2ex}{2.3ex plus .2ex}{\Large\bf}}

\def\maketitle{\par
 \begingroup
 \def\thefootnote{\fnsymbol{footnote}}
 \def\@makefnmark{\hbox
 to 0pt{$^{\@thefnmark}$\hss}}
 \if@twocolumn
 \twocolumn[\@maketitle]
 \else \newpage
 \global\@topnum\z@ \@maketitle \fi\thispagestyle{plain}\@thanks
 \endgroup
 \setcounter{footnote}{0}
 \let\maketitle\relax
 \let\@maketitle\relax
 \gdef\@thanks{}\gdef\@author{}\gdef\@title{}\let\thanks\relax}

\gdef\@title{}
\def\title#1{\gdef\@title{\vskip 25mm
\noindent {\large\bf {#1}\par} \vskip 2em}}

\gdef\@author{}
\def\author#1{\gdef\@author{\noindent
 {\bf \begin{tabular}[t]{c} #1 \end{tabular}\par}
 \vskip 1.5em }}

\gdef\@affil{}
\def\affil#1{\gdef\@affil{\footnotesize\noindent{\em {#1} \par}
 \vskip 1em}}

\gdef\@abstract{}
\def\abstract#1{\gdef\@abstract{\footnotesize
\noindent{\bf Abstract.}{ #1}
\par \vskip 1em}}

\gdef\@keywords{}
\def\keywords#1{\gdef\@keywords{\footnotesize
\noindent{\bf Key Words.}{ #1}
\par \vskip 4em}}

\def\@maketitle{\leftmargini=15mm
\newpage
\null
 \quotation
 \@title
 \@author
 \@affil
 \par
 \vskip 1.5em
 \@abstract
 \@keywords
 \endquotation
 \par
 \vskip 0.5em}

\def\section{\@startsection {section}{1}{\z@}{-3.5ex plus -1ex minus
 -.2ex}{2.3ex plus .2ex}{\normalsize\bf\centering}}
\def\subsection{\@startsection{subsection}{2}{\z@}{-3.25ex plus -1ex minus
 -.2ex}{1.5ex plus .2ex}{\normalsize}}
\def\subsubsection{\@startsection{subsubsection}{3}{\z@}{-3.25ex plus
 -1ex minus -.2ex}{1.5ex plus .2ex}{\small}}

\catcode`@=12
\newtheorem{theorem}{Theorem}[section]
\newtheorem{defi}[theorem]{Definition}
\newtheorem{lemma}[theorem]{Lemma}
\newtheorem{ex}[theorem]{Example}

\def\rit#1{\mbox{${\rm I\hskip -0.2em R}^#1$}}
\def\class{\mbox{$\cal A$}}
\def\norma{\mbox{$\|\cdot\| $}}
\def\classV{\mbox{$\cal V$}}
\def\diag{\mbox{$\rm diag$}}
\def\span{\mbox{$\rm span$}}
\def\qcm{\mbox{$\rm qcm$}}
\def\ovm{\mbox{$\rm ovm$}}
\def\absco{\mbox{$\rm absco$}}
\def\co{\mbox{$\rm co$}}
\def\cl{\mbox{$\rm cl$}}

\def\set1N{\mbox{$\{ 1,2,\ldots ,N\}$}}
\def\NbyN{\mbox{$N\hskip -0.4em\times\hskip -0.4em N$}}

\begin{document}
\author{V.S.\,KOZYAKIN, N.A.\,KUZNETSOV and A.V.\,POKROVSKII\thanks{Current address:
Mathematical Department, University of Queensland, Qld 4072 Australia;
Pokrovskii has been supported by the Australian Research Council Grant
A 89132609.}}

\title{TRANSIENTS IN QUASI--CONTROLLABLE SYSTEMS.\\ OVERSHOOTING, STABILITY
AND INSTABILITY}

\affil{Institute of Information Transmission Problems, Russian
Academy of Sciences, 19 Bolshoi Karetny lane, Moscow 101447,
Russia, e-mail kozyakin@iitp.ru}

\abstract{Families of regimes for control systems are studied possessing the
so called quasi--con\-trol\-\l\-abi\-li\-ty property that is similar to the Kalman
con\-trol\-\l\-abi\-li\-ty property. A new approach is proposed to estimate the degree
of transients overshooting in quasi--controllable systems. This approach is
conceptually related with the principle of bounded regimes absence in the
absolute stability problem. Its essence is in obtaining of constructive
a priori bounds for degree of overshooting in terms of the so called
quasi--con\-trol\-\l\-abi\-li\-ty measure. It is shown that relations between stability,
asymptotic stability and instability for quasi--controllable systems are
similar to those for systems described by linear differential or difference
equations in the case when the leading eigenvalue of the corresponding matrix is
simple. The results are applicable for analysis of transients, classical
absolute stability problem, stability problem for desynchronized systems and
so on.}

\bigskip


\maketitle

\thispagestyle{empty}

{\small

\section{INTRODUCTION}
Stability is recognized as the most important intrinsic attribute of a
real control system. Traditionally, when designing a system one tries
to make it as stable as possible. That is why sometimes one neglect the
fact that a large amount of stability of a system is not a guarantee of
its "good" behavior. The reason is that the stability property
characterizes only the asymptotic behavior of a system and doesn't take
into account system behavior during the so called "transient interval".
As a result, a stable system can have large overshooting or "peaks" in
the transient process that can result in complete failure of a system.
As was noted in \cite{Izm89,MiYo80,OlCi88}
 when the regulator in feedback links is chosen to guarantee
as large a degree of stability as possible then, simultaneously,
overshooting of  system state during the transient process grows, i.e., the
peaking effect is heightened.

Currently, there are a growing number of cases in which systems are described
as operating permanently as if in the transient regime. Examples are flexible
manufacturing systems, adaptive control systems with high level of external
noises, so called desynchronized or asynchronous discrete event systems
\cite{AKKK,BerTsi,KKKK,CHINA} and many others. In connection with this,
it is necessary to develop effective and simple criteria
to estimate the state vector amplitude within the whole time interval of the
system's functioning including the  interval
of transient regime
and an infinite interval when the state vector is "close to equilibrium".

In this paper a new approach is developed presenting the means to solve
effectively the problem of estimation the state vector
amplitude within the whole time interval. The key concept used is a
quasi--con\-trol\-\l\-abi\-li\-ty property of a system that is similar, but weaker than
the Kalman con\-trol\-\l\-abi\-li\-ty property. The
degree of quasi--con\-trol\-\l\-abi\-li\-ty can be characterized by a
numeric value.
The main result of the paper is in proving the following: if a quasi--controllable system
is stable then the amplitudes of all its state trajectories
starting from the unit ball  are bounded by the value reciprocal of
the quasi--con\-trol\-\l\-abi\-li\-ty
measure. Due to the fact that the
measure of quasi--con\-trol\-\l\-abi\-li\-ty can be easily computed, this fact becomes an efficient tool
for analysis
of transients.
Is shown also that for
quasi--controllable systems the properties of stability or instability
are robust with respect to small perturbation of system's parameters.

\section{QUASI-CONTROLLABLE $\cal O$-FAMILIES}
The results of this Section generalize part of results obtained by
Pyatnitskii and Rappoport in \cite{PyatRap}.
Let ${\cal T}$ be either the real half-axes $t\ge 0$ or the set of all
nonnegative integers. The collection ${\cal X}$ is considered, elements of
which are the functions $x(t)$ $(t\in {\cal T})$ taking values in \rit{N}. In
applications functions from ${\cal X}$ often play the role of output signals
for some family of dynamic control systems.

Let us introduce the following notation
$$
{\cal X}(u,t,s)=\{y=x(t):\ x(\cdot )\in {\cal X},\ x(s)=u\},
$$
where $u\in\rit{N}$, $s,t\in{\cal T}$, $s\le t$; a norm in \rit{N} will be denoted by {\norma}.

Suppose that the family ${\cal X}$ possesses the following properties
naturally arising in applications:

\begin{itemize}
\item it satisfies the {\em sewing condition}, i.e. $\forall
\{t_1,t_2,t_3\in{\cal T}: t_1\le t_2\le t_3\}$, $\forall\{x,y\in{\cal X}:
x(t_2)=y(t_2)\}$ $\{\exists z\in{\cal X}: z(t)=x(t)\ {\rm for}\ t_1\le t
\le t_{2},\ z(t)=y(t)\ {\rm for}\ t_2\le t\le t_{3}\}$;

\item it is {\em homogeneous}, i.e., $ {\cal X}(\lambda u,t,0)=\lambda {\cal
X}(u,t,0)$ for $u\in\rit{N}$, $t\in{\cal T}$, $\lambda\in\rit{1}$;

\item it is {\em time invariant}, i.e., ${\cal X}(u,t,s)={\cal X}(u,t-s,0)$
for $u\in\rit{N}$, $s,t\in{\cal T}$, $s\le t$;

\item it is {\em concave}, i.e., ${\cal X}(\lambda u+(1-{\lambda})v,t,0)
\subseteq{\lambda}{\cal X}(u,t,0)+(1-\lambda ){\cal X}(v,t,0)$ for
$u,v\in\rit{N}$, $t\in{\cal T}$, $\lambda\in [0,1]$;

\item it is {\em locally bounded}, i.e., there exists such a continuous
function $\rho (t)$ that the inequality $\|x\|\le{\rho}(t)$ is true for
$x\in{\cal X}(u,t,0)$, $\|u\|=1$, $t\in{\cal T}$;

\item it is {\em continuous}, i.e., from $x_{n}\in\rit{N}$, $x_{n}\rightarrow
x_{0}$ for any $t\in {\cal T}$ follows the equality
$\lim_{n\rightarrow\infty} H\{{\cal X}(x_{n},t,0),{\cal X}(x_{0},t,0)\} =
0$, where $H$ is the Hausdorff distance between sets in \rit{N}.
\end{itemize}

\begin{defi}\label{Ofam} A homogeneous, time invariant, concave, locally bounded and continuous collection of functions satisfying the sewing condition will be
called an $\cal O$-family (output family).
\end{defi}

\begin{ex}\label{EX12} Consider the collection ${\cal X}$ of all solutions of
the class of differential equations of the form
\begin{equation}\label{EqPR}
\frac{dx}{dt} = Ax + bu,\qquad |u|\le\gamma |\langle c,x\rangle |,
\end{equation}
where $b,c\in\rit{N}, c\neq 0$, $A$ is an {\NbyN} matrix and $\gamma $ is a
positive constant. The collection ${\cal X}$ is an {$\cal O$-}family.
\end{ex}

Denote by ${\co}(W)$ and ${\span}(W)$ respectively the convex hull and
the linear span of a vector set $W\subseteq\rit{N}$;
${\absco}(W)={\co}(W\bigcup -W)$ denotes the absolute convex hull of $W$.

\begin{ex}\label{EX13} Let {\class} be a bounded collection of linear
operators acting in \rit{N}. Let us set ${\co}{\class}(x) = {\co}\{\bigcup Ax:\ A\in {\class}\}$, $x\in\rit{N}$. Consider the collection ${\cal X}$ of all solutions of the differential
inclusion
$
{dx}/{dt}\in {\co}{\class}(x).
$
The collection ${\cal X}$ is an $\cal O$-family.
\end{ex}

The last Example covers variety of cases, including Example \ref{EX12},
occurring in the absolute stability theory (see, e.g., \cite{PyatRap}).

\begin{ex}\label{EX14} Let the collection {\class} be the same as in the
previous example. Denote by ${\cal X}$ the collection of piecewise smooth
continuous functions $x(t)$ $(t\ge 0)$ for each of which there exists a
numeric sequence $t_{n}$ $(t_{0}=0,\ t_{n}<t_{n+1},\ n=0,1,2,\ldots )$ such
that the function $x(t)$ satisfies on each interval $[t_{n},t_{n+1}]$ one of
the differential equations of the form
$
{dx}/{dt}=Ax\quad (A\in {\class}).
$
The collection ${\cal X}$ is an $\cal O$-family.
\end{ex}

\begin{ex}\label{EX15} Let ${\cal T}$ be the set of nonnegative integers and
let {\class} be a bounded collection of {\NbyN} matrixes. Let us set
${\class}(x)=\bigcup\{Ax:\ A\in {\class}\}$ for each $x\in\rit{N}$. Consider
the collection ${\cal X}={\cal X}({\class})$ of all solutions of difference
inclusion
$
x(t+1)\in {\class}(x(t)).
$
The collection ${\cal X}$ is an $\cal O$-family.
\end{ex}

This example contains as particular cases a variety of families arising in
the theory of desynchronized  and vertex systems \cite{AKKK,BerTsi,KKKK,CHINA},
two of which
will
be discussed in more details in Examples \ref{EX21}
and \ref{EX211}.

\subsection{\underline{Quasi--con\-trol\-\l\-abi\-li\-ty}} A subspace
$E\subseteq\rit{N}$ is said to be {\em invariant} with respect to the family
${\cal X}$ if for any $x\in E$ the inclusion ${\cal X}(x,t,0)\subseteq E$ is
valid.

\begin{defi}\label{quasicon} An $\cal O$-family ${\cal X}$ will be called
quasi--controllable one if there is no nonzero proper subspace
$E\subset\rit{N}$ invariant with respect to $\cal X$.
\end{defi}

Denote by ${\cal S}$ the set of all $s\in {\cal T}$ such that the amount of
elements of the set $\{t\in {\cal T}:\ 0\le t\le s\}$ is not less than $N$.
Then one can state that the family ${\cal X}$ is quasi--controllable if and
only if for any nonzero $x\in\rit{N}$ and for any $s\in {\cal S}$
\begin{equation}\label{span}
{\span}\{{\cal X}(x,t,0):\ t\in {\cal T},\ 0\le t\le s\} =\rit{N}.
\end{equation}

In the important for the control theory cases of Examples
\ref{EX12}-\ref{EX15}, the quasi--con\-trol\-\l\-abi\-li\-ty property occurs rather often
and it can be verified efficiently. So, in Example \ref{EX12} for
quasi--con\-trol\-\l\-abi\-li\-ty of ${\cal X}$ it is necessary and sufficient that the
pair $\{A,b\}$ be completely controllable and that the pair $\{A,c\}$ be
completely observable by
Kalman. In Examples \ref{EX13}-\ref{EX15} for
quasi--con\-trol\-\l\-abi\-li\-ty of ${\cal X}$ it is necessary and sufficient that the
matrixes from the collection {\class} would not have a common invariant
subspace.

\subsection{\underline{Quasi--con\-trol\-\l\-abi\-li\-ty measure}}
In order to verify the quasi--con\-trol\-\l\-abi\-li\-ty of a family ${\cal X}$ it is
preferable to use in some situations the criteria (\ref{span}) rather than
Definition \ref{quasicon}. However, in a variety of situations the most
suitable tool for verification of quasi--con\-trol\-\l\-abi\-li\-ty of a family ${\cal
X}$ is the numerical measure of quasi--con\-trol\-\l\-abi\-li\-ty proposed below.
Denote the ball of the radius $t$ in a norm {\norma} by ${\bf S}(t)$.

\begin{defi}\label{quasimes} Let $s\in {\cal T}$. The number
$
{\qcm}_{s}({\cal X}) =\inf_{x\in{\rm I\hskip -0.2em R}^{N},\|x\|=1}
\sup\{\rho :{\bf S}(\rho )\subseteq
{\absco}({\cal X}(x,t,0): t\in {\cal T},\ 0\le t\le s)\}
$
will be called the $s$-measure of quasi--con\-trol\-\l\-abi\-li\-ty of the family
${\cal X}$. The number
$$
{\qcm}({\cal X})=\sup_{s\in {\cal S}}{\qcm}_{s}({\cal X})
$$
will be called the quasi--con\-trol\-\l\-abi\-li\-ty measure of the family
${\cal X}$.
\end{defi}

\begin{lemma}\label{Lanalog12} If the family ${\cal X}$ is quasi--controllable,
then ${\qcm}_{s}({\cal X})>0$ for any $s\in {\cal S}$. If ${\qcm}_{s}({\cal
X})>0$ for some $s\in {\cal S}$, then the family ${\cal X}$ is
quasi--controllable.
\end{lemma}

In some situations of interest it is not very difficult to find lower bounds
for quasi--con\-trol\-\l\-abi\-li\-ty measure (see Examples \ref{EX21}, \ref{EX211} below and \cite{KoPo}).

\section{OVERSHOOTING}
Below an efficient and conceptually simple way for estimation of norms of
functions from a given quasi--controllable $\cal O$-family is proposed. The
corresponding estimates are valid within the whole time interval $t\in\cal
T$, including  the initial time interval of transient regime and the
following infinite time interval when the state vector of a system is close
to the equilibrium.
 The principal result may be described as follows: if
quasi--controllable $\cal O$-family  is stable  then
the amplitudes of all its state trajectories starting from the
unit ball  are bounded by the reciprocal of the quasi--con\-trol\-\l\-abi\-li\-ty measure. Due to the fact that the measure
of quasi--con\-trol\-\l\-abi\-li\-ty can be easily computed
in a variety of cases, this criterion is
an efficient tool for analysis of transients in
control systems.

\subsection{\underline{Theorem on an a priori bound}}
\begin{defi}\label{stabLyap}An $\cal O$-family ${\cal X}$ is said to be
Lyapunov stable at the point $x$ if for some $\mu >$0
$$
\sup\{\|u\|:\ u\in\bigcup ({\cal X}(x,t,0),\ t\in {\cal T})\}\le\mu .
$$
An $\cal O$-family ${\cal X}$ is said to be by Lyapunov stable if
for some $\mu >$0
\begin{equation}\label{stabineq}
\sup\{\|u\|:\ u\in\bigcup ({\cal X}(x,t,0),\ t\in {\cal T})\}\le\mu \| x \|.
\end{equation}
An $\cal O$-family ${\cal X}$ is said to be exponentially
Lyapunov stable if for some $\epsilon ,\mu >0$ it is true the inequality
$$
\sup\{\|u\|:\ u\in\bigcup {\cal X}(x,t,0)\}\le{\mu}e^{-\epsilon t}\| x \| .
$$
\end{defi}

The greatest lower bound of those $\mu $ for which the inequality
(\ref{stabineq}) holds is called the {\em overshooting measure} of the
family ${\cal X}$ and is denoted by ${\ovm}({\cal X})$. Note that the
overshooting measure formally characterizes "the value of peaking effect".

\begin{theorem}\label{T21} If the quasi--controllable $\cal O$-family ${\cal X}$ is
stable then
$
{\ovm}({\cal X})\le ({{\qcm}({\cal X})})^{-1}.
$
\end{theorem}

One application of Theorem \ref{T21} to the es\-ti\-ma\-ti\-on of peaking
effect of classical control systems (\ref{EqPR}) immediately follows from
\cite{KoPo}.
Another application is discussed below.

\subsection{\underline{Desynchronized systems}}
Let ${\cal X}={\cal X}({\class})$ be class of $\cal O$-families as in Example
\ref{EX15}, induced by a bounded collection of matrices {\class}. In
desynchronized systems theory (see
\cite{AKKK,KKKK,CHINA}) the collection of matrixes {\class} is often defined as
follows. Given a scalar \NbyN \ matrix $A = (a_{ij})$, then denote $\class = \{A_{1}$,
$A_{2}$, \ldots , $A_{N}\}$ where each $A_i$ is obtaned from the identity matric by replasing $i$-th string by the corresponding string of the matrix $A$.
The matrix $A$ is said to be {\em irreducible}, if by any renumeration of the
basis elements in $\rit{N}$ it cannot be represented in a block triangle form.
%
Let the norm ${\norma}$ in $\rit{N}$ be $\|x\|
=|x_{1}| + |x_{2}| +\ldots + |x_{N}|$. Define the values
$\alpha ={\frac{1}{2N}}\min\{\|(A-I)x\|:\ \|x\|=1\}$,
$\beta ={\frac{1}{2}}\min\{|a_{ij}| :\  i\neq j, a_{ij}\neq 0\}$.

\begin{ex}\label{EX21} The $\cal O$-family ${\cal X}({\class})$ is
quasi--con\-t\-rol\-l\-able if and only if the number $1$ is not an eigenvalue of the
matrix $A$ and $A$ is irreducible. In this case
${\qcm}_{N}[{\cal X}({\class})]\ge{\alpha}{\beta}^{N-1}$.
\end{ex}
\subsection{\underline{Vertex systems}}

Let ${\cal X}={\cal X}({\classV})$ be the $\cal O$-family from Example
\ref{EX15} induced by a bounded collection of matrixes ${\classV}=
\{D_{1}A$, $D_{2}A$, \ldots , $D_{N}A\}$. Here $A$ is an {\NbyN} scalar
matrix with entries $a_{ij}$ and
$ D_{i}= \diag \{d_{1i},\ldots ,d_{ii},\ldots ,d_{Ni}\}$, $i=1,2,\ldots ,N,
$
where $d_{ij}=1$ if $i\neq j$ and $d_{ij}=-1$ if $i=j$.
The $\cal O$-family ${\cal X}({\classV})$ is called  {\em a vertex
family}.

Let again the norm ${\norma}$ in $\rit{N}$ be defined by the equality $\|x\|
=|x_{1}| + |x_{2}| +\ldots + |x_{N}|$. Define the values
$\tilde\alpha ={\frac{1}{N}}\min\{\|Ax\|:\ \|x\|=1\}$,
$\tilde\beta =\min\{|a_{ij}|:\  i\neq j, a_{ij}\neq 0\}$.

\begin{ex}\label{EX211} The {$\cal O$-family} ${\cal X}({\classV})$ is
quasi--con\-trol\-lable if and only if the number $0$ is not an eigenvalue of the
matrix $A$ and matrix $A$ is irreducible. In this case
${\qcm}_{N}[{\cal X}({\classV})]\ge{\tilde\alpha}{\tilde\beta}^{N-1}$.
\end{ex}

\section{ROBUSTNESS} Below it is shown that the quasi--con\-trol\-\l\-abi\-li\-ty property,
as well as properties of stability or instability for quasi--controllable
$\cal O$-families, are robust with respect to small
perturbations.  Due to these properties quasi--controllable $\cal
O$-families are  very attractive object of study for control
theory.

Remind that the symbol $H$ denotes the Hausdorff distance between closed
sets from \rit{N}.
 Say that $\cal O$-families ${\cal X}_{n}$ converge to the $\cal
O$-family ${\cal X}_{*}$ if for all $x\in\rit{N}$, $s\in{\cal T}$ the
following relation is valid
\begin{equation}\label{conv}
\lim_{n\rightarrow\infty }
\sup_{t\le s}H\{{\cl}({\cal X}_{n}(x,t,0)),{\cl}({\cal X}_{*}(x,t,0))\} =
0.
\end{equation}
Let us present a simple auxiliary property of the $s$-measure of
quasi--con\-trol\-\l\-abi\-li\-ty. For any two $\cal O$-families ${\cal X}$ and
${\cal X}_{*}$ let us define the value
$$
\rho _{s}({\cal X},{\cal X}_{*})=\sup_{t\le s,
\| x\|\le 1}H\{{\cl}({\cal X}(x,t,0)),
\ {\cl}({\cal X}_{*}(x,t,0))\},
$$
where $s\in\cal S$. Then
$
|{\qcm}_{s}({\cal X}_{*})-{\qcm}_{s}({\cal X})|\le
\rho _{s}({\cal X},{\cal X}_{*}).
$
This implies the following result:

\begin{theorem}\label{contrrob} Let the sequence of $\cal O$-families ${\cal
X}_{n}$ converge to the quasi--controllable $\cal O$-family ${\cal X}$. Then
${\cal X}_{n}$ is quasi--controllable for sufficiently large $n$
and
$
\lim_{n\rightarrow\infty} {\qcm}_{s}({\cal X}_{n}) = {\qcm}_{s}({\cal X})
$
is valid for each $s\in {\cal S}$.
\end{theorem}

\subsection{\underline{Robustness of stability}}
The next result follows from Theorems \ref{T21} and \ref{contrrob}.
\begin{theorem}\label{T22} Let the sequence of $\cal O$-families ${\cal
X}_{n}$ converge to the quasi--controllable  Lyapunov stable $\cal
O$-family ${\cal X}$. Then
$
\lim_{n\rightarrow\infty} {\ovm}({\cal X}_{n})
\le ({{\qcm}({\cal X})})^{-1},
$
and, for sufficiently large $n$,  ${\cal X}_{n}$ is Lyapunov stable.
\end{theorem}

\begin{theorem}\label{T24} Let the sequence of Lyapunov stable $\cal
O$-families ${\cal X}_{n}$ converges to a quasi--controllable $\cal
O$-family ${\cal X}$. Then the family ${\cal X}$ is also Lyapunov stable.
\end{theorem}

From Theorem \ref{T22} it follows that any $\cal O$-family sufficiently close
to a stable quasi--controllable $\cal O$-family (in the sense
(\ref{conv})) is itself stable and quasi--controllable. At the same time a
stable but not quasi--controllable $\cal O$-family can sometimes be
approximated by stable $\cal O$-families with unbounded quasi--con\-trol\-\l\-abi\-li\-ty
measures. As an example consider the sequence of $\cal O$-families
${\cal X}_{n}$ each of which is the set of solutions of linear difference
vector equations $x(k+1) = A_{n}x(k)$, where $A_{n}$ are two-dimensional
matrixes of the form
$$
A_{n} =\left(\begin{array}{cc}
1-{\epsilon}_{n}&{\delta}_{n}\\0&1-{\epsilon}_{n}
\end{array}\right),
$$
where $\delta _{n}$ tends to zero significantly slower than $\epsilon
_{n}$. The sequence of $\cal O$-families ${\cal X}_{n}$ converges to the
$\cal O$-family ${\cal X}$ consisting of solutions of linear vector equation
$x(k+1) = x(k)$, which is stable but not quasi--controllable.


\subsection{\underline{Robustness of instability}}
An $\cal O$-family ${\cal X}$ will be called {\em absolutely
exponentially unstable with the exponent $\epsilon >0$}, if there exists a
constant $\gamma >0$ such that for each $u\in\rit{N}$, $\| u\| =1$,
there corresponds a function $x(t,u)\in {\cal X}$ satisfying
$
\| x(t,u)\|\ge\gamma e^{\epsilon t}\quad (t\in {\cal T}).
$

The property of absolute exponential instability is a strong one but
unfortunately it is rather difficult to verify. In order for the
quasi--controllable $\cal O$-family ${\cal X}$ to be exponentially unstable,
it is necessary and sufficient that for some $\epsilon >0$ the inequality
$$
\liminf_{t\rightarrow\infty }
\inf_{\|u\|=1}e^{-\epsilon t}\sup\{\|x\|:\ x\in {\cal X}(u,t,0)\} > 0
$$
be fulfilled.

\begin{theorem}\label{T23} Let a quasi--controllable $\cal O$-family ${\cal X}$ be
not Lyapunov stable at least at one point $x$ (see Def.
{\rm\ref{stabLyap}}). Then the family ${\cal X}$ as well as any $\cal O$-family
close to ${\cal X}$ in the sense of convergence {\rm (\ref{conv})} is absolutely
exponentially unstable.
\end{theorem}

\subsection{\underline{Structure of quasi--controllable families}}The previous theorems
can be supplemented by the following two statements on a structure of
quasi--controllable $\cal O$-families.

\begin{theorem}\label{T25} Suppose that each function $x(\cdot)$
from a quasi--controllable
{$\cal O$-}family ${\cal X}$ tends to zero when $t\rightarrow\infty$. Then
the family ${\cal X}$ is exponentially stable.
\end{theorem}

It should be stressed that there are no any assumptions concerning the
closure of the family ${\cal X}$ in the condition of the last Theorem.
The next
classification Theorem follows from Theorems \ref{T25} and \ref{T23}.

\begin{theorem}\label{T26} Every quasi--controllable $\cal O$-family  is either
exponentially stable or is absolutely exponentially unstable, or it is
Lyapunov stable and contains functions $x(t)$ that do not tend to zero for
$t\rightarrow\infty $.
\end{theorem}

Some ideas used in proving the above statements are similar to those used in
proving of the principle of bounded solutions absence \cite{Ko,KrPo}.

\section{ACKNOWLEDGEMENTS} The authors would like to thank Phil Diamond
and Mark Krasnoselskii for
reading a first draft of the manuscript and making many useful suggestions.

}
\end{document}